\documentclass{amsart}
\usepackage{amssymb}
\usepackage{xy}

\newcommand{\GEinfty}{$G$-$E_{\infty}$}
\newcommand{\dL}{\mathrm{L}}
\newcommand{\GSR}{\aS^{G}_{\bR^{\infty}}}
\newcommand{\GSU}{\aS^{G}_{U}}
\newcommand{\Spectra}{\aS}
\newcommand{\nO}{\oE}
\newcommand{\NO}{\oN}
\newcommand{\gO}{\oE_{G}}
\newcommand{\Index}{\aI}
\newcommand{\SQ}{\Xi}

\newcommand{\overto}[1]{\xrightarrow{\,#1\,}}

\newcommand{\putatop}[2]{\genfrac{}{}{0pt}{}{\scriptscriptstyle #1}{\scriptscriptstyle #2}}

\let\Lim\lim
\def\lim{\Lim\nolimits}

\mathchardef\varDelta="7101

\newcommand{\noloc}{\;{:}\,}

\let\iso\cong
\let\sma\wedge

\renewcommand{\to}{\mathchoice{\longrightarrow}{\rightarrow}{\rightarrow}{\rightarrow}}
\newcommand{\from}{\mathchoice{\longleftarrow}{\leftarrow}{\leftarrow}{\leftarrow}}

\DeclareMathAlphabet{\catsymbfont}{U}{rsfs}{m}{n}

\newcommand{\aI}{{\catsymbfont{I}}}

\newcommand{\aS}{{\catsymbfont{S}}}

\newcommand{\oE}{{\mathcal{E}}}
\newcommand{\oN}{{\mathcal{N}}}

\newcommand{\bE}{{\mathbb{E}}}

\newcommand{\bN}{{\mathbb{N}}}
\newcommand{\bS}{{\mathbb{S}}}

\newcommand{\bR}{{\mathbb{R}}}

\def\quickop#1{\expandafter\DeclareMathOperator\csname
#1\endcsname{#1}}
\quickop{id}\quickop{Id}\quickop{holim}\quickop{hocolim}\quickop{op}
\quickop{co}\quickop{Ar}\quickop{sing}\quickop{Hom}\quickop{w}\quickop{Ho}
\quickop{ob}\quickop{diag}\quickop{Cyc}\quickop{Fib}\quickop{Cof}
\quickop{tr}\quickop{trc}\quickop{Gal}\quickop{colim}\quickop{triv}
\quickop{red}\quickop{hoeq}\quickop{inc}\quickop{cyc}
\quickop{fin}\quickop{cy}\quickop{Fix}\quickop{Tot}\quickop{res}\quickop{Sp}
\quickop{Span}\quickop{Gpd}\quickop{Cat}\quickop{Fin}\quickop{Mod}\quickop{CAlg}
\quickop{Vect}\quickop{Fun}\quickop{can}
\quickop{Tr}\quickop{Ad}
\quickop{Aut}\quickop{Map}

\numberwithin{equation}{section}
\newtheorem{thm}[equation]{Theorem}
\newtheorem*{thm*}{Theorem}

\newtheorem{lem}[equation]{Lemma}
\newtheorem{prop}[equation]{Proposition}
\newtheorem{conj}[equation]{Conjecture}
\theoremstyle{definition}
\newtheorem{defn}[equation]{Definition}
\newtheorem{notn}[equation]{Notation}

\xyoption{arrow}
\xyoption{curve}
\xyoption{matrix}
\xyoption{cmtip}
\SelectTips{cm}{}

\newcommand{\term}[1]{\textit{#1}}

\newdir{ >}{{}*!/-5pt/\dir{>}}

\begin{document}

\title%
{A multiplicative version of the tom Dieck splitting}

\author{Andrew J. Blumberg}
\address{Department of Mathematics, Columbia University, 
New York, NY \ 10027}
\email{blumberg@math.columbia.edu}
\author{Michael A. Mandell}
\address{Department of Mathematics, Indiana University,
Bloomington, IN \ 47405}
\email{mmandell@iu.edu}

\ifx\draftdate\undefined\else\date{\draftdate}\fi
\subjclass[2020]{Primary 55P43, 55P91}

\begin{abstract}
While the classical tom Dieck splitting in equivariant stable homotopy
theory is typically regarded as a formula for suspension spectra in
the genuine equivariant stable category, it can be interpreted as a
calculation of the fixed points of $G$-spectra that are derived
pushforwards from the naive equivariant stable category.  We establish a corresponding multiplicative splitting formula for derived
pushforwards to $N_{\infty}$ ring spectra.  Just as the usual tom
Dieck splitting characterizes the equivariant stable category
associated to an $N_{\infty}$ operad $\mathcal{N}$, the multiplicative tom
Dieck splitting characterizes the $G$-symmetric monoidal structure on
the genuine equivariant stable category associated to $\mathcal{N}$.
\end{abstract}

\maketitle

\section{Introduction}


The structure of the equivariant stable category is notoriously
complicated.  Following seminal work of Lewis on incomplete
universes~\cite{Lewis-Splitting}, the modern
understanding of this situation is that the complexity is the price of
having equivariant Poincar\'e duality for smooth $G$-manifolds.
Formally, this can be encoded in a variety of
ways~\cite{LMS,Lewis-Splitting,MM,Blumberg-Thesis,Campion-FreeDuals};
for example, the existence of transfers, the Wirthm\"uller isomorphisms,
invertibility of representation spheres, and Spanier-Whitehead
duality for orbit spectra are all more or less equivalent formulations
of the required structure.

The $G$-equivariant multiplicative structure of genuine equivariant
stable homotopy theory does not currently have such a concise
characterization.  In recent years, the foundational work that went
into the resolution of the Kervaire invariant one problem has brought
new insights into the multiplicative structure of the equivariant
stable category.  In particular, the so-called multiplicative norms
constructed by Hill-Hopkins-Ravenel~\cite{HHR} (building on work of
Greenlees-May~\cite{GreenleesMay-MUCompletion}) appear to characterize
the equivariant multiplicative
structure of the equivariant stable category in a sense that does not
yet admit a precise formulation. For example, recent work on
$N_\infty$ operads~\cite{BlumbergHill-NormsTransfers} connects norms to various kinds of
equivariant commutative ring spectra.  The relationship between
$N_\infty$ operads and $G$-symmetric monoidal structures is more
complicated; see for example~\cite{BlumbergHill-GSymmetric} for preliminary work on
untangling this situation.

In this paper, we revisit a central computational result in genuine
equivariant stable homotopy theory, the tom Dieck splitting.  This
splitting formula describes the categorical fixed points of suspension
spectra, and encodes enough of
the transfers in the case of a finite group to provide another
characterization of the genuine equivariant stable
category among equivariant stable
categories indexed on
universes~\cite[4.14]{BlumbergHill-StableIncomplete}. Moreover, as $H$ ranges
over the subgroups of $G$, the splitting formulas for $H$-fixed points
encode all transfers
and distinguish among even the exotic equivariant stable categories
associated to arbitrary $N_{\infty}$ operads, \textit{ibid}.  In its 
most basic form, the tom Dieck splitting is a computation of the
categorical fixed points of genuine equivariant suspension spectra.
We observe that in fact the tom Dieck splitting holds for all
genuine equivariant spectra in the image of the derived pushforward
functor from the naive equivariant stable category.  
For a precise statement, let $\GSR$ denote the naive equivariant
stable category indexed on the $G$-trivial universe $\bR^{\infty}$ and
let $\GSU$ denote the genuine equivariant stable category on the
complete universe $U$.  Let $i$ denote the inclusion of the trivial
universe in the complete universe and let
\[
i^{\dL}_{*} \colon \Ho(\GSR) \to \Ho(\GSU)
\]
denote the derived change of universe functor.

\begin{thm}\label{main:add}
Let $G$ be a finite group. For $X$ in $\GSR$, there is a natural
isomorphism in the non-equivariant stable category 
\[
\bigvee_{(K)\leq G} (X^{K})_{hWK}\to (i^{\dL}_{*}X)^{G}
\]
where the wedge on the lefthand side is over conjugacy classes of
subgroups of $G$ (choosing one representative in each).
\end{thm}

In the statement, the map $X^{K}\to (i^{\dL}_{*}X)^{G}$ is the
composite
\[
X^{K}\to (i^{\dL}_{*}X)^{K} \to (i^{\dL}_{*}X)^{G}
\]
of the natural map $X^{K}\to (i^{\dL}_{*}X)^{K}$ and the transfer
$(i^{\dL}_{*}X)^{K}\to (i^{\dL}_{*}X)^{G}$.  Because all the functors and natural
transformations commute with desuspension by trivial representations
and homotopy colimits, the previous theorem is an immediate
consequence of the classical tom Dieck splitting on suspension spectra
of orbits.  This theorem (to our knowledge) is due to Lewis and first
appears in~\cite[2.5(b)]{Lewis-Splitting}.  

The main result of this paper constructs a multiplicative analogue of
this generalized tom Dieck splitting. Let $\nO$ denote a non-equivariant 
$E_{\infty}$ operad, which we regard as a $G$-equivariant operad with
trivial $G$-action. Let $\gO$ denote a genuine \GEinfty\
operad. Terminologically, we refer to $\nO$-algebras in $\GSU$ as
$E_{\infty}$ algebras in $\GSU$ and $\gO$-algebras in $\GSU$ as
\GEinfty\ algebras in $\GSU$.  Without loss of generality, we can
assume that we have a map of $G$-equivariant operads
\[
j\colon \nO\to \gO
\]
(representing the unique map in the homotopy category).  
Using brackets to denote algebra categories, $j$ induces a forgetful functor
\[
j^{*}\colon \GSU[\gO]\to \GSU[\nO]
\]
from the category of $\gO$-algebras in $\GSU$ to
the category of $\nO$-algebras in $\GSU$.  This is the right adjoint
of a Quillen adjunction, with left adjoint 
\[
j_{*}\colon \GSU[\nO]\to \GSU[\gO] 
\]
called \term{pushforward}.  The derived pushforward
\[
j^{\dL}_{*}\colon \Ho(\GSU[\nO])\to \Ho(\GSU[\gO])
\]
is the multiplicative analogue of the change of universe
in Theorem~\ref{main:add}: in operadic models for equivariant stable
homotopy theory, change of universe is induced by pushforward.  

In our multiplicative analogue of the generalized tom Dieck splitting,
derived geometric fixed points $X^{\Phi G}$ replace categorical fixed
points $X^{G}$, and when $X=R$ has the structure of a \GEinfty\
algebra or even an $E_{\infty}$ algebra, $R^{\Phi G}$ has the structure
of an $E_{\infty}$ algebra in non-equivariant spectra.  The wedge in
the tom Dieck splitting becomes the derived coproduct of $E_{\infty}$
algebras, which is (up to equivalence) the smash product.  There is a
multiplicative analogue of homotopy orbits using topologically indexed
(homotopy) colimits.  Specifically, we can take the (left derived)
tensor of an $E_{\infty}$ $WK$-spectrum with the $WK$-space $EWK$ and
take the coequalizer in $E_{\infty}$ non-equivariant spectra of the
two $WK$-actions.  The following theorem then gives a multiplicative
tom Dieck splitting.

\begin{thm}\label{main:main}
Let $G$ be a finite group.  Then for $R$ an $E_{\infty}$ algebra in
$\GSU$, there is a natural isomorphism in the homotopy category of 
$E_{\infty}$ algebras in non-equivariant spectra
\[
\coprod_{(K)\leq G} R^{\Phi K}\otimes^{\dL}_{WK}EWK
\to (j^{\dL}_{*}R)^{\Phi G}
\]
where the coproduct on the lefthand side is over conjugacy classes of
subgroups of $G$ (choosing one representative in each).
\end{thm}

This theorem arose from our work on the homotopy theory of the
multiplicative geometric fixed point functor in~\cite{BM-cycl23},
where we stated a multiplicative tom Dieck splitting for \GEinfty\ 
algebras in $\GSU$ that come from non-equivariant $E_{\infty}$
algebras.  The proof of Theorem~\ref{main:main} in
Section~\ref{sec:main} likewise relies on the tools developed
in~\cite{BM-cycl23}, giving a careful study of the more general
case in the statement.  A version of this theorem for the case $G =
C_p$ is implicit in work of Yang~\cite[4.14]{Yang2025Normed}.

Both the additive and multiplicative tom Dieck splittings have
generalizations for pushforwards in the non-genuine setting.  The
statements are the same as the statements of
Theorems~\ref{main:add} and~\ref{main:main} but with the coproduct
indexed over a smaller set of subgroups $K\leq G$.  We state the
multiplicative version in detail. Let  
$\NO$ be any $N_{\infty}$ operad. Without loss of generality, we
assume $\NO$ comes with a map of $G$-equivariant operads
\[
j_{\nO}^{\NO}\colon \nO\to \NO
\]
which induces a Quillen adjunction
\[
\xymatrix@-1pc{%
(j_{\nO}^{\NO})_{*}\colon \GSU[\nO]\ar@<.5ex>[r]&
\ar@<.5ex>[l]\GSU[\NO]\noloc (j_{\nO}^{\NO})^{*}.
}
\]
We then get a derived pushforward $(j_{\nO}^{\NO})^{\dL}_{*}$ from the
homotopy category of $E_{\infty}$ algebras in $\GSU$ to the homotopy
category of $\NO$-algebras in $\GSU$.

\begin{thm}\label{main:incomplete}
Let $G$ be a finite group and $\NO$ an $N_{\infty}$ $G$-operad.  Then
for $R$ an $E_{\infty}$ algebra in $\GSU$, there is a 
natural isomorphism in the homotopy category of 
$E_{\infty}$ algebras in non-equivariant spectra
\[
\coprod_{\putatop{(K)\leq G}{G/K\in \Index_{\NO}}} R^{\Phi K}\otimes^{\dL}_{WK}EWK
\to ((j_{\nO}^{\NO})^{\dL}_{*}R)^{\Phi G}
\]
where the coproduct on the lefthand side is over conjugacy classes of
subgroups of $G$ such that $G/K$
is admissible in the indexing system of $\NO$
(choosing one representative in each).
\end{thm}

We review admissible $G$-sets and indexing systems in Section~\ref{sec:transfer}.

Just as the classical tom Dieck splitting (for $H$-fixed points as $H$
varies over the subgroups of $G$) distinguishes indexing systems and
thereby distinguishes among equivariant stable categories, the
multiplicative tom Dieck splitting in Theorem~\ref{main:incomplete}
(and its instances for $H\leq G$) distinguishes among the possible
$G$-symmetric monoidal structures in genuine $G$-equivariant homotopy
theory.  The idea of a $G$-symmetric monoidal structure on an indexing
system $\Index$ is an enhancement of a symmetric monoidal structure
where the $n$th monoidal power functor generalizes to admissible
$H$-sets in $\Index$.  That is to say, in addition to $n$th powers,
the structure also includes norms.  In particular, the homotopy theory
of commutative monoids in such a structure would be equivalent to the
homotopy theory of $\NO$-algebras, where $\NO$ has the given indexing
system.  It is expected that such structures exist for all indexing
systems $\Index$; they have at least been constructed on the genuine
$G$-equivariant model structure of $G$-equivariant orthogonal spectra
for the indexing systems associated to linear isometries
operads~\cite{BlumbergHill-GSymmetric}.

In the case of compact Lie groups, the situation is significantly more
complicated and less well-understood.  The additive and multiplicative
theories seem qualitatively distinct in that the genuine additive theory
does not admit an operadic model: for a positive dimensional compact
Lie group $G$, categories of $N_{\infty}$ spaces only admit deloopings by 
representations of $\pi_{0}G$ and not by general finite dimensional
representations of $G$ \cite[Ch.~9]{MayMerlingOsorno-memo1540}.
As a consequence, any theory based on operads can only access
transfers associated to subgroups $K\leq H$ of finite index. In contrast,
the additive theory of genuine $G$-spectra admits transfers for
arbitrary closed subgroups, with the transfer of positive codimension
subgroups often coming with a dimension shift.  The classical tom Dieck
splitting for genuine $G$-spectra includes these dimension-shifting
transfer summands.  Multiplicatively, the only known
notion of an equivariant commutative ring spectrum is operadic and we
see only summands for finite index subgroups:

\begin{thm}\label{main:compactLie}
Let $G$ be a compact Lie group.  Let 
\[
j^{\dL}_{*}\colon \Ho(\GSU[\nO])\to \Ho(\GSU[\gO])
\]
be the derived pushforward functor from $E_{\infty}$ algebras in
$\GSU$ to \GEinfty\ algebras in
$\GSU$.  Then for $R$ an $E_{\infty}$ algebra in $\GSU$, there is a
natural isomorphism in the homotopy category of 
$E_{\infty}$ algebras in non-equivariant spectra
\[
\coprod_{\putatop{(K)\leq G}{[G:K]<\infty}} R^{\Phi K}\otimes^{\dL}_{WK}EWK
\to (j^{\dL}_{*}R)^{\Phi G}
\]
where the coproduct on the lefthand side is over conjugacy classes of
finite index subgroups of $G$ (choosing one representative in each).
\end{thm}

We do not currently have a multiplicative theory with positive
codimensional multiplicative transfers.  We discuss this situation in
more detail and make some conjectures about the existence of an
enhanced multiplicative structure on certain commutative ring
$G$-spectra in Section~\ref{sec:speculation}.

\subsection*{Acknowledgments}

This paper reflects ideas coming in part from discussions
with our collaborators Teena Gerhardt, Tyler Lawson, and Mike Hill.
More broadly, the authors are deeply indebted to Gaunce Lewis, Peter
May, and Mike Hopkins for the many things they have taught us over the
years.  The second author would like to thank the Isaac Newton
Institute for Mathematical Sciences, Cambridge, for support and
hospitality during the programme ``Equivariant homotopy theory in
context'' (EHT) where work on this paper was undertaken. This work was
supported by EPSRC grant no.\ EP/R014604/1 and NSF grants 
DMS-2104420, DMS-2104348, DMS-2405029, and DMS-2405030.

\section{Background on $N_{\infty}$ operads, indexing systems, and transfers}
\label{sec:transfer}

We begin with a brief review of indexing systems for $N_{\infty}$
operads and multiplicative transfers for $N_{\infty}$ algebras in
$\GSU$.  Everything in this section is known to experts and it is
provided to extend the definitions in the standard references to the
case of compact Lie groups and finite index subgroups.  In this
section, $G$ denotes a compact Lie group.

Let $H$ be a closed subgroup of $G$.  Given a continuous homomorphism
$\sigma \colon H\to \Sigma_{n}$, we let $H_{\sigma}$ denote the
(closed) subgroup of $G\times \Sigma_{n}$ given by the graph of the
homomorphism $\sigma$:
\[
H_{\sigma}=\{(h,\sigma(h))\mid h\in H\}\leq G\times \Sigma_{n}.
\]
Such a subgroup is called a \term{graph subgroup}. The set of
graph subgroups of $G\times \Sigma_{n}$ is the set of closed subgroups
$\tilde H\leq G\times \Sigma_{n}$ with $\tilde H\cap (1\times
\Sigma_{n})=1$; the homomorphism $\sigma$ is uniquely determined by the
equation $(h,\sigma(h))\in \tilde H$.  We recall that an $N_{\infty}$
$G$-operad is an operad $\NO$ in $G$-spaces such that for each $n$:
\begin{itemize}
\item $\NO(n)^{G}$ is non-empty, 
\item $\NO(n)$ is equivariantly homotopy equivalent
to a $(G\times \Sigma_{n})$-cell complex,
\item the action of the subgroup $1\times \Sigma_{n}\leq G\times
\Sigma_{n}$ on $\NO(n)$ is free, and
\item for every closed subgroup $\tilde H\leq G\times \Sigma_{n}$, the
$\tilde H$-fixed points of $\NO(n)$ are either empty or contractible.
\end{itemize}
It follows that both $\NO$ and $\NO^{G}$ are non-equivariantly $E_{\infty}$
operads, and the fixed points $\NO(n)^{\tilde H}$ are empty when
$\tilde H$ is not a graph subgroup.

We can define admissible $H$-sets in this context just as in the
finite group context, following~\cite[4.5]{BlumbergHill-StableIncomplete}:

\begin{defn}
Let $G$ be a compact Lie group and let $\NO$ be an $N_{\infty}$
$G$-operad.  For $H\leq G$, a finite $H$-set $X$ of cardinality $n$ is \term{admissible} when, given an enumeration of $X$, the graph subgroup
$H_{\sigma }$ associated to the resulting permutation representation
$\sigma \colon H\to\Sigma_{n}$ has a fixed point in $\NO(n)$.
\end{defn}

The definition of indexing system
in~\cite[\S3.2]{BlumbergHill-NormsTransfers}, \cite[1.2]{BlumbergHill-StableIncomplete}
abstracts the properties that the collection of admissible sets
obtains by virtue of the operad structure of~$\NO$. 

\begin{defn}
Let $G$ be a compact Lie group.  An indexing system $\Index$ for $G$
consists of a collection $\Index(H)$ of finite $H$-sets for each
closed subgroup $H\leq G$ such that:
\begin{itemize}
\item $\Index(H)$ contains the one-point set with the trivial action.
\item $\Index(H)$ is closed under isomorphisms, finite disjoint
unions, and finite limits.
\item For each $K\leq H\leq G$ with $K$ finite index in $H$, if $H/K\in
\Index(H)$ and $X\in \Index(K)$ then $H\times_{K} X\in \Index(H)$.
\end{itemize}
\end{defn}

The proof in the finite group case adapts to show that taking
$\Index_{\NO}(H)$ to be the set of admissible $H$-sets for $\NO$
defines an indexing system $\Index_{\NO}$.  As in the finite group
case, it is essentially tautological to show that a map of
$N_{\infty}$ operads is an equivalence (that is, a levelwise
equivariant homotopy equivalence) if and only if they have the same
indexing system, and the cartesian product trick
\[
\NO\from \NO\times \NO'\to \NO'
\]
shows that any pair of operads $\NO$ and $\NO'$ with the same indexing
system are equivalent.  As a consequence, the homotopy category of
$N_\infty$ operads is equivalent to a full subcategory of the poset of
indexing systems.  In the case of finite groups, every indexing system
comes from some $N_{\infty}$ operad~\cite{BonventrePereira,GutierrezWhite,Rubin-Combinatorial}; in the
compact Lie group case, this question has not been studied and it is
not fully established that the techniques in the finite group context
extend. 

The indexing system, or equivalently, the collection of admissible $H$-sets indicates what transfers are inherent in $\NO$-algebras.  To give the most
structure we need to work with an $N_{\infty}$ operad that has a
\term{self-interchange map}: a sequence of maps of $G$-spaces
\[
\NO(m)\times \NO(n)\to \NO(mn)
\]
defining an action of the $G$-operad $\NO$ on itself as in
\cite[\S10]{GuillouMay-EqIterPerm} (the equivariant version of
interchange~\cite[Ch.~2\S3]{BV-Book}).  When $G$ is
finite~\cite[4.18]{Rubin-Categorifying} shows that every equivalence
class of $N_{\infty}$ $G$-operad contains a model with a
self-interchange map; for any compact Lie group, the Steiner operads
admit self-interchange maps.  

The essential point of a self-interchange map is that when $\NO$
admits a self-interchange map, for any $\xi\in \NO(m)$, the action of
$\xi$ on an $\NO$-algebra $A$, 
\[
A\sma \dotsb \sma A\to A
\]
becomes a map of non-equivariant $\NO$-algebras. For $K\leq G$
of codimension zero with $m=[G:K]$, choosing an enumeration of the finite set $G/K$, we
get a permutation
representation $\sigma_{G/K} \colon G\to \Sigma_{m}$ determining the
$G$ action on $G/K$. If $\xi$ is a $G_{\sigma_{G/K}}$-fixed point of
$\NO(m)$, then the map above refines to a $G$-equivariant map of $\NO$-algebras
\[
\mu_{\xi}\colon N_{K}^{G}A\to A.
\]
Note that the refinement depends on the choice of enumeration but is
natural in the algebra $A$.
Taking the point-set geometric $G$-fixed point functor $\Phi^{G}$ (on
the complete universe) as defined in~\cite[\S{}V.4]{MM} and composing
with the HHR diagonal map~\cite[B.209]{HHR}
\[
\Phi^{K}A\to \Phi^{G}N_{K}^{G}A,
\]
we get a map of non-equivariant $\NO^{G}$-algebras
\[
\nu_{\xi}\colon \Phi^{K}A\to \Phi^{G}N_{K}^{G}A\to \Phi^{G}A.
\]
Since (as noted above) when $\NO$ is an $N_{\infty}$ $G$-operad, $\NO^{G}$
is an $E_{\infty}$ operad (in non-equivariant spaces), $\nu_{\xi}$ is
in particular a map of non-equivariant $E_{\infty}$ algebras. 

\begin{defn}[Multiplicative transfer]\label{defn:multtrans}
Let $\NO$ be an $N_{\infty}$ $G$-operad with a self-interchange map,
let $G/K$ be an admissible $G$-set, and let $\sigma_{G/K}\colon G\to
\Sigma_{m}$ be a permutation representation determining the $G$-action
on $G/K$.  For $\xi\in \NO(m)^{G_{\sigma_{G/K}}}$, the \term{multiplicative
$\xi$-transfer} is the natural transformation
\[
\nu_{\xi}\colon \Phi^{K}\to \Phi^{G}
\]
of functors from $\NO$-algebras in $\GSU$ to $\NO^{G}$-algebras in
non-equivariant spectra given as the composite
\[
\Phi^{K}A\to \Phi^{G}N_{K}^{G}A\to \Phi^{G}A
\]
induced by the HHR diagonal and action of $\Phi$ with
$\NO^{G}$-structure induced by the interchange map. The resulting
natural transformation of derived functors 
\[
A^{\Phi K}\to A^{\Phi G}
\]
in the homotopy category of $\NO^{G}$-algebras 
is called the \term{multiplicative transfer for $K\leq G$}.
\end{defn}

Because by definition $\NO(m)^{G_{\sigma_{G/K}}}$ is contractible,
the multiplicative transfer in the homotopy category for $K\leq G$ is
independent of the choice of point $\xi$.  Specifically, the map
$\mu_{\xi}$ is continuous in the element $\xi$, inducing a map 
\[
\NO(m)^{G_{\sigma_{G/K}}}\to \GSU[\NO](N_{K}^{G}A,A)
\]
from the $G_{\sigma_{G/K}}$-fixed points of $\NO(m)$ to the space of
$\NO$-algebra maps in $\GSU$ from $N_{K}^{G}A$ to $A$.  The
point-set geometric fixed point functor is also continuous and so
induces a map
\[
\GSU[\NO](N_{K}^{G}A,A)\to 
\Spectra[\NO^{G}](\Phi^{G}N_{K}^{G}A,\Phi^{G}A)
\to
\Spectra[\NO^{G}](\Phi^{K}A,\Phi^{G}A)
\]
(where $\GSU[\NO]$ denotes the category of $\NO$-algebras in $\GSU$
and  $\Spectra[\NO^{G}]$ denotes the category of
$\NO^{G}$-algebras in non-equivariant spectra).
We have written this out for Theorem~\ref{thm:transequi} below, which
refines the multiplicative transfer to the $WK$-equivariant context.

\begin{thm}\label{thm:transequi}
With notation and assumptions as in Definition~\ref{defn:multtrans},
there exists a natural free $WK$-action on the contractible space
$\NO(m)^{G_{\sigma_{G/K}}}$  such that the composite map
\begin{multline*}
\NO(m)^{G_{\sigma_{G/K}}}\to \GSU[\NO](N_{K}^{G}A,A)\to
\Spectra[\NO^{G}](\Phi^{G}N_{K}^{G}A,\Phi^{G}A)\\
\to
\Spectra[\NO^{G}](\Phi^{K}A,\Phi^{G}A)
\end{multline*}
is $WK$-equivariant for the trivial $WK$-action on $\Phi^{G}A$
and the residual $WK$-action on $\Phi^{K}A$.
\end{thm}

\begin{proof}
Let $\alpha \colon WK\to \Sigma_{m}$ be the permutation representation
for the right action of $WK$ on $G/K$ for the same enumeration used to
construct $\sigma_{G/K}$.  Then essentially by definition, elements of
$\Sigma_{m}$ in the image of $\alpha$ commute with those in the image
of $\sigma_{G/K}$.  It follows that the homomorphism
\[
1\times \alpha \colon WK\to G\times \Sigma_{m}
\]
lands in the normalizer of $G_{\sigma_{G/K}}$.  As a consequence, the
action of $WK$ on $\NO(m)$ via the homomorphism $1\times \alpha$
restricts to the $G_{\sigma_{G/K}}$-fixed points
$\NO(m)^{G_{\sigma_{G/K}}}$; this is the free action in the statement.

As a point-set functor, $N_{K}^{G}A$ is the smash power of $A$ indexed
by $G/K$, with $G$ acting diagonally on $A$ and $G/K$. A
$G$-automorphism of $G/K$ therefore induces an endomorphism of
$N_{K}^{G}A$ as an object of $\GSU$, defining a $WK$-action on
$N_{K}^{G}A$ (in the category $\GSU$).  Because the $\NO$-algebra
structure on $N_{K}^{G}A$ is via the diagonal map $\NO\to \NO^{\times
G/K}$ (which lands in the $WK$-fixed points for the action on the
indexed cartesian power), the $WK$-action on $N_{K}^{G}A$ is through
$\NO$-algebra maps. 

Now let $w\in WK$ and $\xi\in \NO(m)^{G_{\sigma_{G/K}}}$.  We write
$w_{*}$ for the corresponding self-map of $N_{K}^{G}A$ from the
$WK$-action just described. We note that the chosen enumeration
of $G/K$ induces an isomorphism of underlying non-equivariant spectra
$N_{K}^{G}A\iso A^{(m)}$, and under this isomorphism, 
the action of $w$ on $A^{(m)}$ is by the permutation $\alpha(w)$.
Given $\xi \in \NO(m)^{G_{\sigma_{G/K}}}$, the action of $w$ on $\xi$
is by $1\times \alpha(w)$, which is the right operadic permutation
action by $\alpha(w^{-1})$, and so we have 
\[
\mu_{\xi}\circ w^{-1}_{*}=\mu_{\xi \alpha(w^{-1})}=\mu_{w\xi}.
\]
This shows that the map 
\[
\NO(m)^{G_{\sigma_{G/K}}}\to \GSU[\NO](N_{K}^{G}A,A)
\]
is $WK$-equivariant for the actions above with the induced action on
the mapping space where we give the target $A$ the trivial
$WK$-action.  By functoriality, the map
\[
\Phi^{G}\colon \GSU[\NO](N_{K}^{G}A,A)\to
\Spectra[\NO^{G}](\Phi^{G}N_{K}^{G}A,\Phi^{G}A)
\]
is also $WK$-equivariant.  

To complete the proof, it remains to see
that the HHR diagonal map 
\[
\Phi^{K}A\to \Phi^{G}N_{K}^{G}A
\]
is $WK$-equivariant. This does not depend on the algebra structure and
so it suffices to check it for an arbitrary $G$-spectrum in $\GSU$.
By \cite[B.190]{HHR}, we see that the wedge of maps from
$\Phi^{K}F_{V}A(V)$ to $\Phi^{K}A$ is a categorical epimorphism, so it
suffices to consider the case when $A$ is of the form $F_{V}X$.  In
this case, \cite[B.209]{HHR} gives a concrete construction of the map
where the equivariance is clear.
\end{proof}

Since $\NO(m)^{G_{\sigma_{G/K}}}$ is contractible, we get an
essentially unique $WK$-equivariant map $EWK\to
\NO(m)^{G_{\sigma_{G/K}}}$ for any model of $EWK$.  (Indeed, since
$\NO(m)^{G_{\sigma_{G/K}}}$ 
is $WK$-free and non-equivariantly contractible, we can take the map
to be the identity.)  Adjoint to the map in the theorem, we then
have a $WK$-equivariant map of $\NO^{G}$-algebras in $\aS$
\[
\Phi^{K}A\otimes EWK\to \Phi^{G}A,
\]
where $(-)\otimes EWK$ here denotes the point-set indexed colimit with
the space $EWK$.  Since $\Phi^{G}A$ has the trivial $WK$-action, this
factors through the quotient by the $WK$-action in the category of
$\NO^{G}$-algebras in $\aS$,
\begin{equation}\label{eq:hotransfer}
\Phi^{K}A\otimes_{WK} EWK\to \Phi^{G}A.
\end{equation}
We have described a natural point-set map, but taking $A$ to be
cofibrant, $\Phi^{G}$ represents the derived geometric fixed points
$A^{\Phi G}$, and we get a map of derived functors
\begin{equation}\label{eq:derhotrans}
A^{\Phi K}\otimes^{\dL}_{WK} EWK\to A^{\Phi G}.
\end{equation}

For our main theorems, we want to consider the case when $R$ is an
$\NO^{G}$-algebra in $\GSU$ and $A=j_{*}R$ is the pushforward to
$\NO$-algebras in $\GSU$. (Here $j\colon \NO^{G}\to \NO$ is the
inclusion.) The unit of the adjunction $R\to j^{*}j_{*}R$ induces a map of
geometric fixed points $\Phi^{K}R\to \Phi^{K}j_{*}R$, which we can compose
with the map~\eqref{eq:hotransfer} to get a map
\begin{equation}\label{eq:mainpsmap}
\Phi^{K}R\otimes_{WK}EWK\to \Phi^{G}j_{*}R.
\end{equation}
When we take $R$ to be a cofibrant $\NO^{G}$-algebra in $\GSU$,
$j_{*}R$ is a cofibrant $\NO$-algebra in $\GSU$  and represents the
derived pushforward $j^{\dL}_{*}R$.  Composing with the map of derived
functors~\eqref{eq:derhotrans}, we get a natural map of derived functors 
\begin{equation}\label{eq:maindermap}
R^{\Phi K}\otimes^{\dL}_{WK}EWK\to (j^{\dL}_{*}R)^{\Phi G}.
\end{equation}

\section{Tools for geometric fixed points of free $\NO$-algebras}
\label{sec:cycl}

The main step in the proof of the theorems in the introduction is the
calculation of the derived geometric fixed points of free
$\NO$-algebras.  In Chapter~V
of~\cite{BM-cycl23}, we developed tools to do such calculations in the
context of ``generalized orbit desuspension spectra''.  In this
section, we review enough of this theory to cover the special cases we
consider.

In the following statements, we let $\NO$ denote an $N_{\infty}$
operad, let $X$ be an unbased $G$-space, and let $V$ be an orthogonal
$G$-representation (by which we mean a finite dimensional vector space
with linear $G$-action and $G$-invariant inner product).  We use the
notation of~\cite[II.4.6]{MM} and write $F_{V}$ for the functor from
$G$-spaces to orthogonal $G$-spectra left adjoint to the $V$th space
functor. (This functor is denoted as $S^{-V}\sma(-)$ in~\cite{HHR}.)
We write $(-)^{(m)}$ for the $m$th smash power of an orthogonal
$G$-spectrum, with its usual $\Sigma_{m}$-action.  We will write
$X^{m}$ and $V^{m}$ for the cartesian powers of $X$ and $V$ (where we
use orthogonal direct sum for the inner product), with
their usual $\Sigma_{m}$-actions.  We study here the $G$-spectrum
\[
\NO(m)_{+}\sma_{\Sigma_{m}} (F_{V}X_{+})^{(m)}
\iso \NO(m)_{+}\sma_{\Sigma_{m}} F_{V^{m}}X^{m}_{+}.
\]
To interface with the theory in \cite{BM-cycl23}, we note that this is
the generalized orbit desuspension spectrum
(see~\cite[V.1.2]{BM-cycl23}) constructed from the
$(G,\Sigma_{m})$ vector bundle
(see~\cite[V.1.1]{BM-cycl23}) 
\[
\eta \colon \NO(m)\times X^{m}\times V^{m}\to \NO(m)\times X^{m}.
\]
We note that by the hypothesis that $\NO$ is $N_{\infty}$, the action
of $\Sigma_{m}$ on the base space is free.

Theorem~V.2.4 in \cite{BM-cycl23} is the 
basic result to calculate the point-set geometric fixed
point functor $\Phi^{G}$ in this setting. In the special case we are
considering, the statement becomes the following.

\begin{prop}\label{prop:bundlegeofix}
For $\NO$, $V$, and $X$ as above,
\[
\Phi^{G}(\NO(m)_{+}\sma_{\Sigma_{m}}(F_{V}X_{+})^{(m)})\iso 
\biggl(\bigvee_{\sigma \colon G\to \Sigma_{m}} \NO(m)^{G_{\sigma}}_{+}\sma F_{(V^{m})^{G_{\sigma}}}(X^{m})^{G_{\sigma}}_{+}\biggr)/\Sigma_{m},
\]
where the wedge is indexed over the set of continuous homomorphisms
$G\to \Sigma_{m}$ and the superscripts $G_{\sigma}$ denote fixed
points for the corresponding graph subgroup.  In forming the
quotient in the formula, $\Sigma_{m}$ acts on the smash powers,
cartesian powers, and $\NO(m)$ in the usual way and acts on the set of
homomorphisms by conjugation.
\end{prop}

The isomorphism in the previous proposition is natural in $X$ 
and $V$. Moreover, the isomorphism does not depend on the fact that
the $\Sigma_{m}$-free $(G\times \Sigma_{m})$-space $\NO(m)$ came from
an operad: it is  natural in equivariant maps of the $\Sigma_{m}$-free
$(G\times \Sigma_{m})$-space $\NO(m)$. 

In the case when $\NO=\NO^{G}$, we have that $\NO(m)^{G_{\sigma}}$ is
empty unless $\sigma$ is the trivial homomorphism and the formula
above simplifies to 
\begin{multline*}
\Phi^{G}(\NO^{G}(m)_{+}\sma_{\Sigma_{m}}(F_{V}X_{+})^{(m)})\iso 
\NO^{G}(m)_{+}\sma_{\Sigma_{m}}F_{(V^{m})^{G}}(X^{m})^{G}_{+}\\
\iso 
\NO^{G}(m)_{+}\sma_{\Sigma_{m}}(F_{V^{G}}X^{G}_{+})^{(m)}.
\end{multline*}
We get the corresponding formula when we take geometric $K$-fixed
points and naturality gives the identification $WK$-equivariantly.
Wedging over $m$, we then get the following observation.

\begin{prop}\label{prop:einffree}
Let $\NO$, $V$, and $X$ be as above, and let $\nO=\NO^{G}$.  Let $\bE$
denote the free $\nO$-algebra functor in either $\GSU$ or $\Spectra$.
Then 
\[
\Phi^{K}\bE(F_{V}X_{+})\iso \bE(F_{V^{K}}X^{K}_{+}),
\]
and this isomorphism is compatible with the residual $WK$-actions on
both sides.
\end{prop}

The formula is not quite as nice for $\NO$ a general $N_{\infty}$
operad, but we can reorganize the
righthand side of Proposition~\ref{prop:bundlegeofix} in
terms of admissible $G$-sets.  Given a continuous homomorphism $\sigma
\colon G\to \Sigma_{m}$, let $S_{\sigma}$ denote the $G$-set obtained
by giving $\{1,\dotsc,m\}$ the corresponding action.  Then (by
definition) $S_{\sigma}$ is admissible exactly when the
$G_{\sigma}$-fixed points of $\NO(m)$ are non-empty.  Moreover,
$S_{\sigma}$ and $S_{\sigma'}$ are isomorphic $G$-sets exactly when
the homomorphisms $\sigma$ and $\sigma'$ are conjugate by the
$\Sigma_{m}$-action in Proposition~\ref{prop:bundlegeofix}, and the group of
$G$-automorphisms of $S_{\sigma}$ is exactly the isotropy subgroup of
$\sigma$ for this conjugation action.  Letting
$C_{m}$ be a set consisting of a single choice of representative for
each conjugacy class, we can then rewrite the righthand side of the
formula in Proposition~\ref{prop:bundlegeofix} as
\begin{equation}\label{eq:simp1}
\bigvee_{\sigma \in C_{m}} \bigl(\NO(m)^{G_{\sigma}}_{+}\sma
F_{(V^{m})^{G_{\sigma}}}(X^{m})^{G_{\sigma}}_{+}\bigr)/\Aut_{G}(S_{\sigma}).
\end{equation}
Since conjugacy classes of homomorphisms correspond to isomorphism
classes of $G$-sets, we can organize this in terms of $G$-sets
instead.  We do so using the following notation.

\begin{notn}\label{notn:vector}
Let $k$ be the number of conjugacy classes of subgroups $K\leq G$ with
$G/K$ admissible.  Choose one subgroup in each conjugacy class and
enumerate them
\[
K_{1},\dotsc,K_{k}.
\]
Write $\vec n=(n_{1},\dotsc,n_{k})$ for an ordered list of $k$ natural numbers (including
zero).  Define the cardinality $\#\vec n$ of $\vec n$ by
\[
\#\vec n=\sum_{i=1}^{k}[G:K_{i}]n_{i},
\]
the $G$-set $S_{\vec n}$ corresponding to $\vec n$ as
\[
S_{\vec n}=(G/K_{1}\times n_{1})\amalg \dotsb \amalg (G/K_{k}\times n_{k}),
\]
and the permutation representation $\sigma_{\vec n}\colon G\to
\Sigma_{\#\vec n}$ corresponding to $\vec n$ as 
\[
\sigma_{\vec n}=(\underbrace{\sigma_{G/K_{1}}\oplus \dotsb \oplus \sigma_{G/K_{1}}}_{n_{1}}) 
\oplus \dotsb \oplus (\underbrace{\sigma_{G/K_{k}}\oplus \dotsb \oplus \sigma_{G/K_{k}}}_{n_{k}}),
\]
where $\oplus$ denotes block sum, and $\sigma_{G/K_{i}}$ is the
permutation representation associated to $G/K_{i}$ chosen in the
previous section.  For $i=1,\dotsc,k$, let $\vec s_{i}$ be the list
with $1$ in the $i$th spot and $0$ elsewhere, so that 
\[
\#\vec s_{i}=[G:K_{i}], \quad S_{\vec s_{i}}=G/K_{i}, 
\quad \sigma_{\vec s_{i}}=\sigma_{G/K_{i}}.
\]
\end{notn}

By construction $\sigma_{\vec n}$ is the permutation representation of
$S_{\vec n}$ for the evident enumeration of $S_{\vec n}$.  The set
$\{S_{\vec n}\}$ consists of one choice of $G$-set in each isomorphism
class, and the set $\{\sigma_{\vec n}\mid \#\vec n=m\}$ gives a choice
for $C_{m}$ above, the set of conjugacy class representatives of
homomorphisms $G\to \Sigma_{m}$.  We can then rewrite the
expression~\eqref{eq:simp1} as 
\begin{equation}\label{eq:simp2}
\bigvee_{\#\vec n=m} \bigl(\NO(m)^{G_{\sigma_{\vec n}}}_{+}\sma
F_{(V^{m})^{G_{\sigma_{\vec n}}}}(X^{m})^{G_{\sigma_{\vec
n}}}_{+}\bigr)/\Aut_{G}(S_{\sigma_{\vec n}}).
\end{equation}
This can be written more intrinsically in terms of $S_{\vec n}$ using
the relationship between the $G_{\sigma_{\vec n}}$-fixed points of a
cartesian product and the space of $G$-maps out of $S_{\vec n}$:
\[
(X^{m})^{G_{\sigma_{\vec n}}}\iso \Map_{G}(S_{\vec n},X)\iso 
(X^{K_{1}})^{n_{1}}\times \dotsb \times (X^{K_{k}})^{n_{k}}.
\]
Analogous observations apply to $(V^{m})^{G_{\sigma_{\vec n}}}$, and
we have a canonical isomorphism 
\begin{multline*}
F_{(V^{m})^{G_{\sigma_{\vec n}}}}(X^{m})^{G_{\sigma_{\vec n}}}_{+}
\iso 
F_{(V^{K_{1}})^{n_{1}}\oplus \dotsb \oplus (V^{K_{k}})^{n_{k}}}
\bigl((X^{K_{1}})^{n_{1}}\times \dotsb \times (X^{K_{k}})^{n_{k}}\bigr)_{+}
\\
\iso
(F_{V^{K_{1}}}X^{K_{1}}_{+})^{(n_{1})}
  \sma \dotsb \sma (F_{V^{K_{k}}}X^{K_{k}}_{+})^{(n_{k})}.
\end{multline*}
The group of $G$-automorphisms of $S_{\vec n}$ is the product of the
groups of $G$-automorphisms of $G/K_{i}\times {n_{i}}$, which 
can be written as a semidirect product of the symmetric group
$\Sigma_{n_{i}}$ with $G$-automorphisms of $G/K_{i}$:
\[
\Aut_{G}(S_{\vec n})\iso 
(\Sigma_{n_{1}}\wr WK_{1})
  \times \dotsb \times 
(\Sigma_{n_{k}}\wr WK_{k}).
\]
The action of $\Aut_{G}(S_{\vec n})$ on $F_{(V^{m})^{G_{\sigma_{\vec
n}}}}(X^{m})^{G_{\sigma_{\vec n}}}_{+}$ corresponds to the evident
action of the symmetric groups and Weyl groups on the smash product
under the isomorphism above.  This then gives the following
up-to-homotopy simplification of Proposition~\ref{prop:bundlegeofix}
and expression~\eqref{eq:simp2}. 

\begin{prop}\label{prop:simp}
Let $\NO$, $V$, and $X$ be as above, and let $\nO=\NO^{G}$.  
Let $\bN$ denote the free $\NO$-algebra functor in $\GSU$ and let
$\bE$ denote the free $\nO$-algebra functor in $\Spectra$.
Then 
\[
\Phi^{G}(\bN (F_{V}X_{+}))\iso \bigvee_{\vec n}\biggl(
\NO(\#\vec n)^{G_{\sigma_{\vec n}}}_{+}
\mathbin{\mathop\sma\limits_{\Aut_{G}(S_{\vec n})}}
\bigl(
(F_{V^{K_{1}}}X^{K_{1}}_{+})^{(n_{1})}
  \sma \dotsb \sma (F_{V^{K_{k}}}X^{K_{k}}_{+})^{(n_{k})}
\bigr)\biggr);
\]
moreover, the map of $\nO$-algebras
\[
\bE \biggl(\bigvee_{i=1}^{k} 
\bigl(\NO([G:K_{i}])^{G_{\sigma_{G/K_{i}}}}_{+}
  \sma_{WK_{i}} F_{V^{K_{i}}}X^{K_{i}}_{+}\bigr)
\biggr)
\to
\Phi^{G}(\bN (F_{V}X_{+}))
\]
induced by the inclusion of the summands $\vec n=\vec s_{i}$
is a homotopy equivalence of the underlying spectra.
\end{prop}

\begin{proof}
We prove the last statement by observing that the map induced by
multiplication in~$\NO$ 
\begin{multline*}
\NO(n_{1}+\dotsb+n_{k})^{G}\times 
(\NO([G:K_{1}])^{G_{\sigma_{G/K_{1}}}})^{n_{1}}
\times \dotsb \times
(\NO([G:K_{k}])^{G_{\sigma_{G/K_{k}}}})^{n_{k}}
\\
\to \NO(\#\vec n)^{G_{\sigma_{\vec n}}}
\end{multline*}
is a homotopy equivalence of free $\Aut_{G}(S_{\vec n})$-spaces. Since
the underlying spaces on both sides are contractible, it is enough to
observe that they are both free.  The space on the righthand side is
free because it is a subspace of $\NO(\#\vec n)$ and the action of
$\Aut_{G}(S_{\vec n})$ is as a subgroup of $\Sigma_{\#\vec n}$.  The
space on the lefthand side is a subspace of
\[
\NO(n_{1}+\dotsb+n_{k})\times 
\NO([G:K_{1}])^{n_{1}}
\times \dotsb \times
\NO([G:K_{k}])^{n_{k}}
\]
which
\[
(\Sigma_{n_{1}}\wr \Sigma_{[G:K_{1}]})\times \dotsb \times (\Sigma_{n_{k}}\wr \Sigma_{[G:K_{k}]})
\]
acts freely on.
\end{proof}

So far we have discussed the point-set geometric fixed point functor.
Theorem~V.2.7 of \cite{BM-cycl23} gives a criterion for
the point-set functor to represent the derived functor.  Our example
of $\NO(m)_{+}\sma_{\Sigma_{m}} (F_{V}X_{+})^{(m)}$ fits this
criterion because for the $(G,\Sigma_{m})$ vector bundle $\eta$ above, the
$(G,\Sigma_{m})$ vector bundle $\eta(G)$ in the statement has base
space
\[
\coprod_{\sigma \colon G\to \Sigma_{m}} \NO(m)^{G_{\sigma}}\times (X^{m})^{G_{\sigma}},
\]
which $\Sigma_{m}$ acts freely on.  The statement in this case is the following.

\begin{prop}\label{prop:dergeo}
The map from the derived geometric fixed points to the point-set
geometric fixed points
\[
(\NO(m)_{+}\sma_{\Sigma_{m}}(F_{V}X_+)^{(m)})^{\Phi G}
\to 
\Phi^{G}(\NO(m)_{+}\sma_{\Sigma_{m}}(F_{V}X_+)^{(m)})
\]
is an isomorphism in the non-equivariant stable category.
\end{prop}

In the main case of interest, $X$ will be a $G$-CW complex and then the
$G$-spectrum 
$\NO(m)_{+}\sma_{\Sigma_{m}}(F_{V}X_+)^{(m)}$ is homotopy equivalent to
a cofibrant $G$-spectrum in the usual
model structure for $G$-equivariant orthogonal spectra on the complete
universe, and Proposition~\ref{prop:dergeo} is a special case
of~\cite[V.4.17]{MM}. 

\section{Proof of main theorems}
\label{sec:main}

In this section, we prove Theorems~\ref{main:main},
\ref{main:incomplete}, and~\ref{main:compactLie} from the
introduction.  Noting that Theorem~\ref{main:main} is also the special case
of Theorem~\ref{main:incomplete} for \GEinfty\ operads and is the
special case of Theorem~\ref{main:compactLie} for zero-dimensional Lie
groups, we can prove all three theorems at once by proving a
theorem generalizing both Theorem~\ref{main:incomplete} and
Theorem~\ref{main:compactLie}.  We can state this theorem more
precisely here, having constructed the basic maps underlying it in
Section~\ref{sec:transfer}.

\begin{thm}\label{main:omnibus}
Let $G$ be a compact Lie group, $\NO$ an $N_{\infty}$ $G$-operad, and
let $\nO=\NO^{G}$.  Let $R$ be an $\nO$-algebra in $\GSU$ and let 
$A=j^{\dL}_{*}R$ be the derived pushforward $\NO$-algebra in $\GSU$.
Indexing on one choice of subgroup $K\leq G$ in each conjugacy class in
which $G/K$ is admissible, the coproduct of the natural
maps~\eqref{eq:maindermap} in the homotopy category of non-equivariant
$\nO$-algebras 
\[
\coprod
R^{\Phi K}\otimes^{\dL}_{WK}EWK\to A^{\Phi G}
\]
is an isomorphism.
\end{thm}

We chose and enumerated the conjugacy class representatives in
Notation~\ref{notn:vector}.  The following then gives a technical
restatement of the previous theorem closer to the tools we apply from
Section~\ref{sec:cycl}. 

\begin{thm}\label{thm:maintech}
With notation as in Theorem~\ref{main:omnibus}, assume that $R$ is
cofibrant as an $\nO$-algebra in $\GSU$.  The
maps~\eqref{eq:mainpsmap} in the point-set 
category of non-equivariant $\nO$-algebras
\[
\Phi^{K_{i}}R\otimes_{WK_{i}}EWK_{i}\to \Phi^{G}(j_{*}R), \qquad i=1,\dotsc,k
\]
induce an isomorphism on the derived coproduct
\[
\coprod_{i=1}^{k}\mathstrut^{\dL}\ \Phi^{K_{i}}R\otimes_{WK_{i}}EWK_{i}
\to
\coprod_{i=1}^{k}\ \Phi^{K_{i}}R\otimes_{WK_{i}}EWK_{i}
\to \Phi^{G}(j_{*}R).
\]
\end{thm}

Theorem~\ref{thm:maintech} implies Theorem~\ref{main:omnibus} because
under the cofibrancy hypothesis, $j_{*}R$ represents the derived pushforward
$j_{*}^{\dL}R$, the geometric fixed points $\Phi^{K_{i}}R$ and
$\Phi^{G}j_{*}R$ represent the derived geometric fixed points $R^{\Phi
K_{i}}$ and $(j_{*}R)^{\Phi G}$, and the indexed colimit $(-)\otimes_{WK_{i}}EWK_{i}$
represents the derived indexed colimit $(-)\otimes^{\dL}_{WK_{i}}EWK_{i}$.

As in the previous section, let $\bN$ denote the free $\NO$-algebra
functor in $\GSU$ and let $\bE$ denote the free $\nO$-algebra
functor in $\GSU$ or $\Spectra$.  Because cofibrant $\nO$-algebras in
$\GSU$ are retracts of cell algebras built from cells of the form
\[
\bE(F_{V}(G/H\times \partial D^{q})_{+})\to 
\bE(F_{V}(G/H\times D^{q})_{+}),
\]
an easy homotopy colimit argument reduces the general case of
Theorem~\ref{thm:maintech} to the case when $R=\bE(F_{V} X_{+})$ for
$X=G/H\times \partial D^{q}$ or $X=G/H\times D^{q}$, and $V$ is an
orthogonal $G$-representation.  In this case,
according to Proposition~\ref{prop:einffree}, the geometric fixed
points simplify, and we have
\[
\Phi^{K_{i}}\bE(F_{V}X_{+})\iso \bE(F_{V^{K_{i}}}X^{K_{i}}_{+})
\]
and
\[
\Phi^{K_{i}}\bE(F_{V}X_{+})\otimes_{WK_{i}}EWK_{i}
\iso \bE(F_{V^{K_{i}}}(X^{K_{i}}\times_{WK_{i}}EWK_{i})_{+}).
\]
Moreover, since these are free $\nO$-algebras in $\Spectra$, the
derived coproduct is easy to calculate, and we have
\begin{multline*}
\coprod\mathstrut^{\dL}\ \Phi^{K_{i}}\bE(F_{V}X_{+})\otimes_{WK_{i}}EWK_{i}
\overto{\iso}
\coprod\mathstrut\ \Phi^{K_{i}}\bE(F_{V}X_{+})\otimes_{WK_{i}}EWK_{i}\\
\iso
\bE\biggl(\bigvee F_{V^{K_{i}}}(X^{K_{i}}\times_{WK_{i}}EWK_{i})_{+}\biggr).
\end{multline*}
Thus, it
suffices to prove the following lemma.

\begin{lem}\label{lem:maintech}
For an unbased $G$-space $X$ and an orthogonal $G$-representation $V$, the map
\begin{multline*}
\bE\biggl(\bigvee_{i=1}^{k} F_{V^{K_{i}}}(X^{K_{i}}\times_{WK_{i}}EWK_{i})_{+}\biggr)
\iso
\coprod\mathstrut\ \Phi^{K_{i}}\bE(F_{V}X_{+})\otimes_{WK_{i}}EWK_{i}\\
\to \Phi^{G}(\bN F_{V}X_{+})
\end{multline*}
induced by the maps~\eqref{eq:mainpsmap} is a weak equivalence.
\end{lem}

The statement in Lemma~\ref{lem:maintech} looks a lot like the
statement in Proposition~\ref{prop:simp}, particularly when we model
the map $EWK\to \NO(m)^{G_{\sigma_{G/K}}}$ by the identity map
in~\eqref{eq:hotransfer}.  To prove Lemma~\ref{lem:maintech}, we need
to identify the map in Proposition~\ref{prop:simp} as the map induced
by~\eqref{eq:hotransfer}.  Because both maps are maps of
$\nO$-algebras in $\Spectra$, the maps
\[
\bE \biggl(\bigvee_{i=1}^{k} 
\bigl(\NO([G:K_{i}])^{G_{\sigma_{G/K_{i}}}}_{+}\sma_{WK_{i}} F_{V^{K_{i}}}X^{K_{i}}_{+}\bigr)
\biggr)
\to
\Phi^{G}(\bN (F_{V}X_{+}))
\]
are determined by the composite maps of spectra
\begin{multline*}
\NO([G:K_{i}])^{G_{\sigma_{G/K_{i}}}}_{+}\sma_{WK_{i}}
F_{V^{K_{i}}}X^{K_{i}}_{+}
\\\to
\bE \biggl(\bigvee_{i=1}^{k} 
\bigl(\NO([G:K_{i}])^{G_{\sigma_{G/K_{i}}}}_{+}\sma_{WK_{i}} F_{V^{K_{i}}}X^{K_{i}}_{+}\bigr)
\biggr)
\to \Phi^{G}(\bN (F_{V}X_{+})),
\end{multline*}
and on the point-set level, these maps are determined by the
restriction 
\begin{multline}\label{eq:restriction}
F_{V^{K_{i}}}X^{K_{i}}_{+}\iso
\{\xi\}_{+}\sma 
F_{V^{K_{i}}}X^{K_{i}}_{+}\\
\to
\NO([G:K_{i}])^{G_{\sigma_{G/K_{i}}}}_{+}\sma_{WK_{i}}
F_{V^{K_{i}}}X^{K_{i}}_{+}
\to \Phi^{G}(\bN (F_{V}X_{+})),
\end{multline}
for each $\xi\in \NO([G:K_{i}])^{G_{\sigma_{G/K_{i}}}}$.  The following
proposition now completes the proof.

\begin{prop}
Let $K=K_{i}$ for some $i=1,\dotsc,k$, and let $m=[G:K]$.  For any $\xi\in
\NO(m)^{G_{\sigma_{G/K}}}$, the restriction~\eqref{eq:restriction} 
\[
F_{V^{K}}X^{K}_{+}\to \Phi^{G}(\bN (F_{V}X_{+}))
\] 
of the map in Proposition~\ref{prop:simp} is the composite
of the map
\[
F_{V^{K}}X^{K}_{+}\to \Phi^{K}(\bE (F_{V}X_{+}))
\]
induced by Proposition~\ref{prop:einffree}, the natural map
\[
\Phi^{K}(\bE (F_{V}X_{+}))\to \Phi^{K}(\bN (F_{V}X_{+})), 
\]
and the map 
\[
\nu_{\xi}\colon \Phi^{K}(\bN (F_{V}X_{+}))\to  \Phi^{G}(\bN (F_{V}X_{+}))
\]
of Definition~\ref{defn:multtrans}.
\end{prop}

\begin{proof}
In terms of the computation in Proposition~\ref{prop:bundlegeofix} of
$\Phi^{G}$ of the summands of $\bN F_{V}X_{+}$, the
map~\eqref{eq:restriction} for Proposition~\ref{prop:simp} sends
$F_{V^{K}}X^{K}_{+}$ into 
\[
\NO(m)^{G_{\sigma_{G/K}}}_{+}\sma
  F_{(V^{m})^{G_{\sigma_{G/K}}}}(X^{m})^{G_{\sigma_{G/K}}}_{+}
\]
by the map induced by the inclusion of $\xi$ and the isomorphisms 
\[
(V^{m})^{G_{\sigma_{G/K}}}\iso V^{K},\qquad
(X^{m})^{G_{\sigma_{G/K}}}\iso X^{K}.
\]

On the other hand, the map $\nu_{\xi}$ is induced by the HHR
diagonal and the geometric $G$-fixed points of the map 
\[
N_{K}^{G}(\bN F_{V}X_{+})\to \bN F_{V}X_{+}
\]
induced by $\xi$ for the $\NO(m)$-multiplication.  Precomposing with
the map obtained by applying $N_{K}^{G}$ to the inclusion of $F_{V}X$
in $\bN F_{V}X$, we can identify the composite map 
\begin{equation}\label{eq:mtcomp}
N_{K}^{G}(F_{V}X_{+})\to
N_{K}^{G}(\bN F_{V}X_{+})\to \bN F_{V}X_{+}
\end{equation}
as follows.  Write $O\xi$ for the $(G\times \Sigma_{m})$-orbit of $\xi$
in $\NO(m)$.  Using the distinguished point $\xi$, we have a canonical
isomorphism $O\xi\iso (G\times \Sigma_{m})/G_{\sigma_{G/K}}$, inducing
an isomorphism 
\[
N_{K}^{G}(F_{V}X_{+})\iso O\xi_{+}\sma_{\Sigma_{m}}(F_{V}X_{+})^{(m)}
\]
(see, for example, \cite[VI.7.1ff]{BM-cycl23}).  In
this notation, \eqref{eq:mtcomp} is the map
\[
N_{K}^{G} F_{V}X_{+}\iso
O\xi_{+}\sma_{\Sigma_{m}} (F_{V}X_{+})^{(m)}
\to \NO(m)_{+}\sma_{\Sigma_{m}}(F_{V}X_{+})^{(m)}
\to \bN F_{V}X_{+}
\]
induced by the inclusion of $\xi$ in $\NO(m)$.  Taking geometric
$G$-fixed points, written in terms of the formula of
Proposition~\ref{prop:bundlegeofix}, we get the map
\begin{multline*}
\Phi^{G}(N_{K}^{G}(F_{V}X_{+}))\iso
\biggl(\bigvee_{\sigma \colon G\to \Sigma_{m}} O\xi^{G_{\sigma}}_{+}\sma F_{(V^{m})^{G_{\sigma}}}(X^{m})^{G_{\sigma}}_{+}\biggr)/\Sigma_{m}\\
\to 
\biggl(\bigvee_{\sigma \colon G\to \Sigma_{m}} \NO(m)^{G_{\sigma}}_{+}\sma F_{(V^{m})^{G_{\sigma}}}(X^{m})^{G_{\sigma}}_{+}\biggr)/\Sigma_{m}
\to \Phi^{G}(\bN F_{V}X_{+}).
\end{multline*}
The HHR diagonal is the (isomorphic) inclusion of $F_{V^{K}}X^{K}_{+}$
as the summand corresponding to $\sigma=\sigma_{G/K}$.  This is the
map described in the previous paragraph for~\eqref{eq:restriction} in
the context of Proposition~\ref{prop:simp}.
\end{proof}

\section{Speculations on structures beyond operads}\label{sec:speculation}

The additive and multiplicative structures on the $G$-stable category
in the setting of finite groups are essentially controlled by the same
formal structures; transfers become norms and the additive tom Dieck
splitting is replaced by a multiplicative tom Dieck splitting of the
same basic form.  In the case of positive dimensional compact Lie
groups, however, additive and multiplicative structures no longer
have parallel theories.

We can explain precisely the nature of the difference using work of
Blumberg-Hill~\cite{BlumbergHill-StableIncomplete}. This work constructs
models of the equivariant stable category with additive structure
controlled by arbitrary $N_\infty$ operads, for $G$ a finite group.
For the Steiner operad on a universe $U$ (the equivariant analogue of
the $E_{\infty}$ little cubes operad), the theory recovers the
classical theory of equivariant $G$-spectra indexed on the universe
$U$, but other $N_{\infty}$ operads produce exotic equivariant stable
categories not equivalent to the classical ones.  The equivariant
stable category for the $N_{\infty}$ operad $\gO$ comes with a
stabilization (suspension spectrum) functor $\Sigma^{\infty}_{\gO}$
on based $G$-spaces, generalizing the classical stabilization
functor. One of the main theorems of the paper is a version of the tom
Dieck splitting for $\Sigma^{\infty}_{\gO}X$.

Everything in~\cite{BlumbergHill-StableIncomplete} works just as well
to construct an equivariant stable category associated to an $N_{\infty}$
operad in the context of a positive dimensional compact Lie group.
For the \GEinfty\ operad $\gO$, the associated $G$-equivariant stable
category is not the genuine $G$-equivariant stable category, but
rather the $G$-equivariant stable category on the universe $U_{0}$
of $\pi_{0}G$-representations.  We get the following tom Dieck
splitting: 
\[
(\Sigma^{\infty}_{\gO}X)^{G}\simeq 
\bigvee_{\putatop{(K)\leq G}{[G:K]<\infty}} \Sigma^{\infty}(X^{K})_{hWK}.
\]
On the other hand, the tom Dieck splitting in the genuine equivariant stable
category has dimension-shifting transfer summands:
\[
(\Sigma^{\infty}_{U}X)^{G}\simeq 
\bigvee_{(K)\leq G} \Sigma^{\infty}(\Sigma^{\Ad_{WK}} X^{K})_{hWK},
\]
where $\Ad_{WK}$ is the $WK$-representation given by the tangent space
of $WK$ at the identity (see, for example, \cite[XIX.1.3]{May-Alaska}).

The multiplicative tom Dieck splitting for $\gO$ above
(Theorem~\ref{main:compactLie}) only has finite
index subgroup summands.
This phenomenon raises the question of whether there exists a multiplicative 
structure analogous to the additive structure of the genuine
equivariant stable category in the case of positive dimensional
compact Lie groups.  Such a structure by necessity would go
beyond what is accessible in the theory of operadic algebras.  For
purposes of discussion, we call this putative structure ``structure
$\SQ$''.  We take the remaining 
paragraphs of the paper to discuss some ideas; these remarks
have benefited from discussions of the authors with Mike Hill.

We can work out what the structure might look like in the special case
of suspension spectra on free infinite loop spaces.  Let $X$ be a
$G$-connected based $G$-space and consider the free genuine
infinite loop space on $X$ (relative to the base point)
\[
Q_{U} X=\Omega^{\infty}_{U}\Sigma^{\infty}_{U}X=\bigcup_V \Omega^{V}\Sigma^{V}X.
\]
Let
$A=\Sigma^{\infty}_{U}(Q_{U}X)_{+}$.  This is at least a \GEinfty\
algebra in $\GSU$, but it also admits some kind of multiplicative
dimension-shifting transfer on geometric fixed points, constructed as
follows.  For $K \leq G$, we have (non-equivariantly)
\[
A^{\Phi K}\simeq \Sigma^{\infty}((Q_{U}X)^{K})_{+}\simeq\Sigma^{\infty}(\Omega^{\infty}(\Sigma^{\infty}_{U}X)^{K})_{+}
\]
and the transfer 
\[
(\Sigma_{U}^{\infty}X)^{K}\sma S^{\Ad_{WK}}\to (\Sigma^{\infty}_{U}X)^{G}
\]
induces a map of non-equivariant $E_{\infty}$ ring spectra
\[
\Sigma^{\infty}(\Omega^{\infty}((\Sigma_{U}^{\infty}X)^{K}\sma S^{\Ad_{WK}}))_{+}
\to
\Sigma^{\infty}(\Omega^{\infty}(\Sigma^{\infty}_{U}X)^{G})_{+}\simeq A^{\Phi G}.
\]
We would like to interpret the object on the left as a suspension (based tensor)
of $A^{\Phi K}$ with $S^{\Ad_{WK}}$ in the category of augmented
$E_{\infty}$ ring spectra, $A^{\Phi K}\mathbin{\widehat\otimes}^{\dL}
S^{\Ad_{WK}}$.  Making this interpretation precise seems to involve questions related
to those that arise in equivariant infinite loop space theory for
compact Lie groups, as studied in~\cite{Blumberg-Thesis}.

Nonetheless, this example suggests the structure in the following
conjecture.

\begin{conj}
There exists a multiplicative structure $\Xi$ so that $\Xi$-algebras
are in particular \GEinfty\ algebras and in the augmented
or nonunital setting admit (dimension-shifting) multiplicative
transfers for all closed subgroups.  For an augmented $E_{\infty}$
algebra $R$, the multiplicative transfers induce a multiplicative tom
Dieck splitting of the form
\[
\coprod_{(K)\leq G} (R^{\Phi K} 
\mathbin{\widehat\otimes}^{\dL} S^{\Ad_{WK}})\otimes^{\dL}_{WK}EWK \overto{\sim}
   ((j_{\nO}^{\Xi})^{\dL}_{*}R)^{\Phi G}.
\]
\end{conj}

From the discussion above and the form dimension-shifting transfers
take additively, it seems best to start the search for structure $\SQ$
in the augmented or non-unital context.  Current work on genuine equivariant
factorization homology for \GEinfty\ algebras in $\GSU$ suggests 
norm / geometric fixed point formulas of the form 
\[
(N_{e}^{G}A)^{\Phi G}\simeq \bS
\]
\cite[2.5]{BHM-Norms} whereas for suspension spectra, the
Wirthm\"uller isomorphism  suggests formulas of the form
\[
(N_{e}^{G}A)^{\Phi G}\simeq A\mathbin{\widehat\otimes} S^{\Ad_{G}},
\]
which resembles equivariant nonabelian
Poincar\'e duality or a multiplicative version of the Wirthm\"uller
isomorphism.

Non-equivariantly, the cobordism spectra provided the prototypes for
developing the notion of $E_{\infty}$ algebras. If we look at the
cobordism spectra equivariantly for compact Lie groups, we see further
hints of a structure beyond the operadic: $BO$ and $BU$ are the zeroth
spaces of genuine equivariant spectra; if structure $\SQ$ exists, we
should expect to see it on equivariant $MO$ and $MU$, and more
generally on Thom spectra associated to infinite loop maps to $BO$.
We can already see such a structure in the augmented setting, looking
for example at $MU\sma MU$ as an augmented $MU$-algebra, using the
Thom isomorphism $MU\sma MU\iso MU\sma BU_{+}$ and the space-level
structure on $BU$.


\bibliographystyle{plain}
\bibliography{bluman}

\end{document}